# CONSISTENCY OF SPECTRAL CLUSTERING

By Ulrike von Luxburg, Mikhail Belkin and Olivier Bousquet

*Max Planck Institute for Biological Cybernetics, Ohio State University
and Pertinence*

Consistency is a key property of all statistical procedures analyzing randomly sampled data. Surprisingly, despite decades of work, little is known about consistency of most clustering algorithms. In this paper we investigate consistency of the popular family of spectral clustering algorithms, which clusters the data with the help of eigenvectors of graph Laplacian matrices. We develop new methods to establish that, for increasing sample size, those eigenvectors converge to the eigenvectors of certain limit operators. As a result, we can prove that one of the two major classes of spectral clustering (normalized clustering) converges under very general conditions, while the other (unnormalized clustering) is only consistent under strong additional assumptions, which are not always satisfied in real data. We conclude that our analysis provides strong evidence for the superiority of normalized spectral clustering.

**1. Introduction.** Clustering is a popular technique which is widely used in statistics, computer science and various data analysis applications. Given a set of data points, the goal is to separate the points in several groups based on some notion of similarity. Very often it is a natural mathematical model to assume that the data points have been drawn from an underlying probability distribution. In this setting it is desirable that clustering algorithms should satisfy certain basic consistency requirements:

- In the large sample limit, do the clusterings constructed by the given algorithm "converge" to a clustering of the whole underlying space?
- If the clusterings do converge, is the limit clustering a reasonable partition of the whole underlying space, and what are the properties of this limit clustering?

Interestingly, while extensive literature exists on clustering and partitioning (e.g., see Jain, Murty and Flynn [27] for a review), very few clustering algorithms have been analyzed or shown to converge in the setting where the data is sampled from an arbitrary probability distribution. In a parametric setting, clusters are often identified with the individual components of a mixture distribution. Then clustering reduces to standard parameter estimation, and of course there exist many results on the consistency of such estimators. However, in a nonparametric setting there are only two major classes of clustering algorithms where convergence questions have been studied at all: single linkage and $k$-means.

The $k$-means algorithm clusters a given set of points in $\mathbb{R}^d$ by constructing $k$ cluster centers such that the sum of squared distances of all data points to their closest cluster centers is minimized (e.g., Section 14.3 of Hastie, Tibshirani and Friedman [23]). Pollard [38] shows consistency of the global minimizer of the objective function for $k$-means clustering. However, as the $k$-means objective function is highly nonconvex, the problem of finding its global minimum is often infeasible. As a consequence, the guarantees on the consistency of the minimizer are purely theoretical and do not apply to existing algorithms, which use local optimization techniques. The same problem also concerns all the follow-up articles on Pollard [38] by many different authors.

Linkage algorithms construct a hierarchical clustering of a set of data points by starting with each point being a cluster, and then successively merging the two closest clusters (e.g., Section 14.3 of Hastie, Tibshirani and Friedman [23]). For this class of algorithms, Hartigan [22] demonstrates a weaker notion of consistency. He proves that the algorithm will identify certain high-density regions, but he does not obtain a general consistency result.

In our opinion, the results about the consistency of clustering algorithms which can be found in the literature are far from satisfactory. This lack of consistency guarantees is especially striking as clustering algorithms are widely used in most scientific disciplines which deal with data in any form.

In this paper we investigate the limit behavior of the class of spectral clustering algorithms. Spectral clustering is a popular technique going back to Donath and Hoffman [17] and Fiedler [19]. In its simplest form, it uses the second eigenvector of the graph Laplacian matrix constructed from the affinity graph between the sample points to obtain a partition of the samples into two groups. Different versions of spectral clustering have been used for many different problems such as load balancing (Van Driessche and Roose [46]), parallel computations (Hendrickson and Leland [24]), VLSI design (Hagen and Kahng [21]) and sparse matrix partitioning (Pothen, Simon and Liou [40]). Laplacian-based clustering algorithms also have found success in applications to image segmentation (Shi and Malik [43]), text mining (Dhillon [15]) and as general purpose methods for data analysis and clustering (Alpert [2], Kannan, Vempala and Vetta [28], Ding et al. [16], Ng, Jordan and Weiss [36] and Belkin and Niyogi [10]). A nice survey on the history of spectral clustering can be found in Spielman and Teng [44]; for a tutorial introduction to spectral clustering, see von Luxburg [48].

While there has been some theoretical work on properties of spectral clustering on finite point sets (e.g., Spielman and Teng [44], Gauttery and Miller [20], Kannan, Vempala and Vetta [28]), so far there have not been any results discussing the limit behavior of spectral clustering for samples drawn from some underlying probability distribution. In the current article, we establish consistency results and convergence rates for several versions of spectral clustering. To prove those results, the main step is to establish the convergence of eigenvalues and eigenvectors of random graph Laplace matrices for growing sample size. Interestingly, our analysis shows that while one type of spectral clustering ("normalized") is consistent under very general conditions, another popular version of spectral clustering ("unnormalized") is only consistent under some very specific conditions which do not have to be satisfied in practice. We therefore conclude that the normalized clustering algorithm should be the preferred method in practical applications.

From a mathematical point of view, the question of convergence of spectral clustering boils down to the question of convergence of spectral properties of random graph Laplacian matrices constructed from sample points. The convergence of eigenvalues and eigenvectors of certain random matrices has already been studied extensively in the statistics community, especially for random matrices of fixed size such as sample covariance matrices, or for random matrices with i.i.d. entries (see Bai [6] for a review). However, those results cannot be applied in our setting, as the size of the graph Laplacian matrices ($n \times n$) increases with the sample size $n$, and the entries of the random graph Laplacians are not independent from each other. In the machine learning community, the spectral convergence of positive definite "kernel matrices" has attracted some attention, as can be seen in Shawe-Taylor et al. [42], Bengio et al. [12] and Williams and Seeger [50]. Here, several authors build on the work of Baker [7], who studies numerical solutions of integral equations by deterministic discretizations of integral operators on a grid. However, his methods cannot be carried over to our case, where integral operators are discretized by a random selection of sample points (see Section II.10 of von Luxburg [47] for details). Finally, Koltchinskii [30] and Koltchinskii and Giné [31] have obtained convergence results for random discretizations of integral operators which are close to what we would need for spectral clustering. However, to apply their techniques and results, it is necessary that the operators under consideration are Hilbert–Schmidt, which turns out not to be the case for the unnormalized Laplacian. Consequently, to prove consistency results for spectral clustering, we have to derive new methods which hold under more general conditions than all the methods mentioned above. As a by-product we recover certain results from Koltchinskii [30] and Koltchinskii and Giné [31] by using considerably simpler techniques.

There has been some debate on the question whether normalized or unnormalized spectral clustering should be used. Recent papers using the normalized version include Van Driessche and Roose [46], Shi and Malik [43], Kannan, Vempala and Vetta [28], Ng, Jordan and Weiss [36] and Meila and Shi [33], while Barnard, Pothen and Simon [8] and Gauttery and Miller [20] use unnormalized clustering. Comparing the empirical behavior of both approaches, Van Driessche and Roose [46] and Weiss [49] find some evidence that the normalized version should be preferred. On the other hand, under certain conditions, Higham and Kibble [25] advocate for the unnormalized version. It seems difficult to resolve this question from purely graph-theoretic considerations, as both normalized and unnormalized spectral clustering can be justified by similar graph theoretic principles (see next section). In our work we now obtain the first theoretical results on this question. They show the superiority of normalized spectral clustering over unnormalized spectral clustering from a statistical point of view.

This paper is organized as follows: In Section 2 we briefly introduce the family of spectral clustering algorithms, and describe what the difference between "normalized" and "unnormalized" spectral clustering is. After giving an informal overview of our consistency results in Section 3, we introduce mathematical prerequisites and notation in Section 4. The convergence of normalized spectral clustering is stated and proved in Section 5, and rates of convergence are proved in Section 6. In Section 7 we establish conditions for the convergence of unnormalized spectral clustering. Those conditions are studied in detail in Section 8. In particular, we investigate the spectral properties of the limit operators corresponding to normalized and unnormalized spectral clustering, point out some important differences, and show theoretical and practical examples where the convergence conditions in the unnormalized case are violated.

**2. Spectral clustering.**    The purpose of this section is to briefly introduce the class of spectral clustering algorithms. For a comprehensive introduction to spectral clustering and its various derivations, explanations and properties, we refer to von Luxburg [48]. Readers who are familiar with spectral clustering or who first want to get an overview over our results are encouraged to jump to Section 3 immediately.

Assume we are given a set of data points $X_1, \ldots, X_n$ and pairwise similarities $k_{ij} := k(X_i, X_j)$ which are symmetric (i.e., $k_{ij} = k_{ji}$) and nonnegative. We denote the data similarity matrix as $K := (k_{ij})_{i,j=1,\ldots,n}$ and define the matrix $D$ to be the diagonal matrix with entries $d_i := \sum_{j=1}^n k_{ij}$. Spectral clustering uses matrices which have been studied extensively in spectral graph theory, so-called graph Laplacians. Graph Laplacians exist in three different flavors. The *unnormalized graph Laplacian* (sometimes also called the combinatorial Laplacian) is defined as the matrix

$$L = D - K.$$

The *normalized graph Laplacians* are defined as

$$L' = D^{-1/2}LD^{-1/2} = I - D^{-1/2}KD^{-1/2},$$

$$L'' = D^{-1}L = I - D^{-1}K.$$

Given a vector $f = (f_1, \ldots, f_n)^t \in \mathbb{R}^n$, the following key identity can be easily verified:

$$f^t L f = \tfrac{1}{2} \sum_{i,j=1}^{n} w_{ij}(f_i - f_j)^2.$$

This equation shows that $L$ is positive semi-definite. It can easily be seen that the smallest eigenvalue of $L$ is 0, and the corresponding eigenvector is the constant one vector $\mathbb{1} = (1, \ldots, 1)^t$. Similar properties can be shown for $L'$ and $L''$. There is a tight relationship between the spectra of the two normalized graph Laplacians: $v$ is an eigenvector of $L''$ with eigenvalue $\lambda$ if and only if $w = D^{1/2}v$ is an eigenvector of $L'$ with eigenvalue $\lambda$. So from a spectral point of view, the two normalized graph Laplacians are equivalent. A discussion of various other properties of graph Laplacians can be found in the literature; see, for example, Chung [14] for the normalized and Mohar [35] for the unnormalized case.

There exist two major versions of spectral clustering, which we call "normalized" or "unnormalized" spectral clustering, respectively. The basic versions of those algorithms can be summarized as follows:

---

**Basic spectral bi-clustering algorithms**
**Input: Similarity matrix** $K \in \mathbb{R}^{n \times n}$.
  Find the eigenvector $v$ corresponding to the second smallest eigenvalue for one of the following problems:

$$Lv = \lambda v \ \text{(for unnormalized spectral clustering)},$$

$$L''v = \lambda v \ \text{(for normalized spectral clustering)}.$$

**Output: Clusters** $A = \{j; v_j \geq 0\}$ and $\bar{A} = \{j; v_j < 0\}$.

---

It is not straight forward to see why the clusters produced by those algorithms are useful in any way. The roots of spectral clustering lie in spectral graph theory. Here we consider the "similarity graph" induced by the data, namely, the graph with adjacency matrix $K$. On this graph, clustering reduces to the problem of graph partitioning: we want to find a partition of the graph such that the edges between different groups have very low weights (which means that points in different clusters are dissimilar from each other) and the edges within a group have high weights (which means that points within the same cluster are similar to each other). Different ways of formulating and solving the objective functions of such graph partitioning problems lead to normalized and unnormalized spectral clustering, respectively. For details, we refer to von Luxburg [48].

Note that the spectral clustering algorithms as presented above are simplified versions of spectral clustering. The implementations used in practice can differ in various details. In particular, in the case when one is interested in obtaining more than two clusters, one typically uses not only the second but also the next few eigenvectors to construct a partition. Moreover, more complicated rules can be used to construct a partition from the coordinates of the eigenvectors than simply thresholding the eigenvector at 0. For details, see von Luxburg [48].

**3. Informal statement of our results.** In this section we want to present our main results in a slightly informal but intuitive manner. For the precise mathematical details and proofs, we refer to the following sections. The goal of this article is to study the behavior of normalized and unnormalized spectral clustering on random samples when the sample size $n$ is growing. In Section 2 we have seen that spectral clustering partitions a given sample $X_1, \ldots, X_n$ according to the coordinates of the first eigenvectors of the (normalized or unnormalized) Laplace matrix. To stress that the Laplace matrices depend on the sample size $n$, from now on we denote the unnormalized and normalized graph Laplacians by $L_n$, $L'_n$ and $L''_n$ (instead of $L$, $L'$ and $L''$ as in the last section). To investigate whether the various spectral clustering algorithms converge, we will have to establish conditions under which the eigenvectors of the Laplace matrices "converge." To see which kind of convergence results we aim at, consider the case of the second eigenvector $(v_1, \ldots, v_n)^t$ of $L_n$. It can be interpreted as a function $f_n$ on the discrete space $\mathcal{X}_n := \{X_1, \ldots, X_n\}$ by defining $f_n(X_i) := v_i$, and clustering is then performed according to whether $f_n$ is smaller or larger than a certain threshold. It is clear that in the limit for $n \to \infty$, we would like this discrete function $f_n$ to converge to a function $f$ on the whole data space $\mathcal{X}$ such that we can use the values of this function to partition the data space. In our case it will turn out that this space can be chosen as $C(\mathcal{X})$, the space of continuous functions on $\mathcal{X}$. In particular, we will construct a degree function $d \in C(\mathcal{X})$ which will be the "limit" of the discrete degree vector $(d_1, \ldots, d_n)$. Moreover, we will explicitly construct linear operators $U$, $U'$ and $U''$ on $C(\mathcal{X})$ which will be the limit of the discrete operators $L_n$, $L'_n$ and $L''_n$. Certain eigenvectors of the discrete operators are then proved to "converge" (in a certain sense to be explained later) to eigenfunctions of those limit operators. Those eigenfunctions will then be used to construct a partition of the whole data space $\mathcal{X}$.

In the case of normalized spectral clustering it will turn out that this limit process behaves very nicely. We can prove that, under certain mild conditions, the partitions constructed on finite samples converge to a sensible partition of the whole data space. In meta-language, this result can be stated as follows:

RESULT 1 (*Convergence of normalized spectral clustering*).    Under mild assumptions, if the first $r$ eigenvalues $\lambda_1, \ldots, \lambda_r$ of the limit operator $U'$ satisfy $\lambda_i \neq 1$ and have multiplicity 1, then the same holds for the first $r$ eigenvalues of $L'_n$ for sufficiently large $n$. In this case, the first $r$ eigenvalues of $L'_n$ converge to the first $r$ eigenvalues of $U'$ a.s., and the corresponding eigenvectors converge a.s. The clusterings constructed by normalized spectral clustering from the first $r$ eigenvectors on finite samples converge almost surely to a limit clustering of the whole data space.

In the unnormalized case, the convergence theorem looks quite similar, but there are some subtle differences that will turn out to be important.

RESULT 2 (*Convergence of unnormalized spectral clustering*).    Under mild assumptions, if the first $r$ eigenvalues of the limit operator $U$ have multiplicity 1 and do not lie in the range of the degree function $d$, then the same holds for the first $r$ eigenvalues of $\frac{1}{n}L_n$ for sufficiently large $n$. In this case, the first $r$ eigenvalues of $\frac{1}{n}L_n$ converge to the first $r$ eigenvalues of $U$ a.s., and the corresponding eigenvectors converge a.s. The clusterings constructed by unnormalized spectral clustering from the first $r$ eigenvectors on finite samples converge almost surely to a limit clustering of the whole data space.

On the first glance, both results look very similar: if first eigenvalues are "nice," then spectral clustering converges. However, the difference between Results 1 and 2 is what it means for an eigenvalue to be "nice." For the convergence statements to hold, in Result 1 we only need the condition $\lambda_i \neq 1$, while in Result 2 the condition $\lambda_i \notin \mathrm{rg}(d)$ has to be satisfied. Both conditions are needed to ensure that the eigenvalue $\lambda_i$ is isolated in the spectrum of the limit operator, which is a fundamental requirement for applying perturbation theory to the convergence of eigenvectors. We will see that in the normalized case, the limit operator $U'$ has the form $Id - T$ where $T$ is a compact linear operator. As a consequence, the spectrum of $U'$ is very benign, and all eigenvalues $\lambda \neq 1$ are isolated and have finite multiplicity. In the unnormalized case, however, the limit operator will have the form $U = M - S$, where $M$ is a multiplication operator and $S$ a compact integral operator. The spectrum of $U$ is not as nice as the one of $U'$, and, in particular, it contains the continuous interval $\mathrm{rg}(d)$. Eigenvalues of this operator will only be isolated in the spectrum if they satisfy the condition $\lambda \notin \mathrm{rg}(d)$. As the following result shows, this condition has important consequences.

RESULT 3 [*The condition $\lambda \notin \mathrm{rg}(d)$ is necessary*].

1. There exist examples of similarity functions such that there exists no nonzero eigenvalue outside of $\mathrm{rg}(d)$.
2. If this is the case, the sequence of second eigenvalues of $\frac{1}{n}L_n$ computed by any numerical eigensolver converges to $\min d(x)$. The corresponding eigenvectors do not yield a sensible clustering of the data space.
3. For a large class of reasonable similarity functions, there are only finitely many eigenvalues (say, $r_0$) inside the interval $]0, \min d(x)[$. In this case, the same problems as above occur if the number $r$ of eigenvalues used for clustering satisfies $r > r_0$.
4. The condition $\lambda \notin \mathrm{rg}(d)$ refers to the limit case and, hence, cannot be verified on the finite sample.

This result complements Result 2. The main message is that there are many examples where the conditions of Result 2 are not satisfied, that in this case unnormalized spectral clustering fails completely, and that we cannot detect on a finite sample whether the convergence conditions are satisfied or not.

To further investigate the statistical properties of normalized spectral clustering, we compute rates of convergence. Informally, our result is the following:

RESULT 4 (*Rates of convergence*). The rates of convergence of normalized spectral clustering can be expressed in terms of regularity conditions of the similarity function $k$. For example, for the case of the widely used Gaussian similarity function $k(x, y) = \exp(-\|x - y\|^2/\sigma^2)$ on $\mathbb{R}^d$, we obtain a rate of $\mathcal{O}(1/\sqrt{n})$.

Finally, we show how our theoretical results influence the results of spectral clustering in practice. In particular, we demonstrate differences between the behavior of normalized and unnormalized spectral clustering.

**4. Prerequisites and notation.** In the rest of the paper we always make the following general assumptions:

GENERAL ASSUMPTIONS. The data space $\mathcal{X}$ is a compact metric space, $\mathcal{B}$ the Borel $\sigma$-algebra on $\mathcal{X}$, and $P$ a probability measure on $(\mathcal{X}, \mathcal{B})$. Without loss of generality we assume that the support of $P$ coincides with $\mathcal{X}$. The sample points $(X_i)_{i \in \mathbb{N}}$ are drawn independently according to $P$. The similarity function $k : \mathcal{X} \times \mathcal{X} \to \mathbb{R}$ is supposed to be symmetric, continuous and bounded away from 0 by a positive constant, that is, there exists a constant $l > 0$ such that $k(x, y) > l$ for all $x, y \in \mathcal{X}$.

Most of those assumptions are standard in the spectral clustering literature. We need the symmetry of the similarity function in order to be able to represent our data by an *undirected* graph (note that spectral graph theory does not carry over to directed graphs as, e.g., the graph Laplacians are no longer symmetric). The continuity of $k$ is needed for robustness reasons: small changes in the data should not change the result too much. For the same reason, we make the assumption that $k$ should be bounded away from 0. This becomes necessary when we consider normalized graph Laplacians, where we divide by the degree function and still want the result to be robust with respect to small changes in the underlying data. Only the compactness of $\mathcal{X}$ is added for mathematical convenience. Most results in this article are also true without compactness, but their proofs would require a serious technical overhead which does not add to the general understanding of the problem.

For a finite sample $X_1, \ldots, X_n$, which has been drawn i.i.d. according to $P$, and a given similarity function $k$ as in the General assumptions, we denote the similarity matrix by $K_n = (k(X_i, X_j))_{i, j \leq n}$ and the degree matrix $D_n$ as the diagonal matrix with the degrees $d_i = \sum_{j=1}^n k(X_i, X_j)$ on the diagonal. Similarly, we will denote the unnormalized and normalized Laplace matrices by $L_n = D_n - K_n$ and $L'_n = D_n^{-1/2} L_n D_n^{-1/2}$. The eigenvalues of the Laplace matrices $0 = \lambda_1 \leq \lambda_2 \leq \cdots \leq \lambda_n$ will always be ordered in increasing order, respecting multiplicities. The term "first eigenvalue" always refers to the trivial eigenvalue $\lambda_1 = 0$. Note that throughout the whole paper, we will use superscript-$t$ (such as $f^t$) to denote the transpose of a vector or a matrix, while "primes" (as in $L'$ or $L''$) are used to distinguish different matrices and operators. $I$ is used to denote the identity matrix.

For a real-valued function $f$, we always denote the range of the function by $\mathrm{rg}(f)$. If $\mathcal{X}$ is connected and $f$ is continuous, $\mathrm{rg}(f) = [\inf_x f(x), \sup_x f(x)]$. The restriction operator $\rho_n : C(\mathcal{X}) \to \mathbb{R}^n$ denotes the (random) operator which maps a function to its values on the first $n$ data points, that is, $\rho_n f = (f(X_1), \ldots, f(X_n))^t$.

Now we want to recall certain facts from spectral and perturbation theory. For more details, we refer to Chatelin [13], Anselone [3] and Kato [29]. By $\sigma(T) \subset \mathbb{C}$, we denote the spectrum of a bounded linear operator $T$ on some Banach space $E$. We define the *discrete spectrum* $\sigma_d$ to be the part of $\sigma(T)$ which consists of all isolated eigenvalues with finite algebraic multiplicity, and the *essential spectrum* $\sigma_{ess}(T) = \sigma(T) \setminus \sigma_d(T)$. The essential spectrum is always closed, and the discrete spectrum can only have accumulation points on the boundary to the essential spectrum. It is well known (e.g., Theorem IV.5.35 in Kato [29]) that compact perturbations do not affect the essential spectrum, that is, for a bounded operator $T$ and a compact operator $V$, we have $\sigma_{ess}(T + V) = \sigma_{ess}(T)$. A subset $\tau \subset \sigma(T)$ is called isolated if there exists an open neighborhood $M \subset \mathbb{C}$ of $\tau$ such that $\sigma(T) \cap M = \tau$. For an isolated part $\tau \subset \sigma(T)$, the corresponding spectral projection $\Pr_\tau$ is defined as $\frac{1}{2\pi i} \int_\Gamma (T - \lambda I)^{-1} \, d\lambda$, where $\Gamma$ is a closed Jordan curve in the complex plane separating $\tau$ from the rest of the spectrum. In the special case where $\tau = \{\lambda\}$ for an isolated eigenvalue $\lambda$, $\Pr_\tau$ is a projection on the invariant subspace related to $\lambda$. If $\lambda$ is a simple eigenvalue (i.e., it has algebraic multiplicity 1), then the spectral projection $\Pr_\tau$ is a projection on the eigenfunction corresponding to $\lambda$.

DEFINITION 5 (*Convergence of operators*).   Let $(E, \| \cdot \|_E)$ be an arbitrary Banach space, $B$ its unit ball, and $(S_n)_n$ a sequence of bounded linear operators on $E$:

- $(S_n)_n$ *converges pointwise*, denoted by $S_n \xrightarrow{p} S$, if $\|S_n x - Sx\|_E \to 0$ for all $x \in E$.
- $(S_n)_n$ *converges compactly*, denoted by $S_n \xrightarrow{c} S$, if it converges pointwise and if for every sequence $(x_n)_n$ in $B$, the sequence $(S - S_n)x_n$ is relatively compact (has compact closure) in $(E, \| \cdot \|_E)$.
- $(S_n)_n$ *converges in operator norm*, denoted by $S_n \xrightarrow{\|\cdot\|} S$, if $\|S_n - S\| \to 0$, where $\| \cdot \|$ denotes the operator norm.
- $(S_n)_n$ is called *collectively compact* if the set $\bigcup_n S_n B$ is relatively compact in $(E, \| \cdot \|_E)$.
- $(S_n)_n$ *converges collectively compactly*, denoted by $S_n \xrightarrow{cc} S$, if it converges pointwise and if there exists some $N \in \mathbb{N}$ such that the operators $(S_n - S)_{n > N}$ are collectively compact.

Both operator norm convergence and collectively compact convergence imply compact convergence. The latter is enough to ensure the convergence of spectral properties in the following sense (cf. Proposition 3.18 and Sections 3.6 and 5.1 in Chatelin [13]):

PROPOSITION 6 (*Perturbation results for compact convergence*).   *Let $(E, \| \cdot \|_E)$ be an arbitrary Banach space and $(T_n)_n$ and $T$ bounded linear operators on $E$ with $T_n \xrightarrow{c} T$. Let $\lambda \in \sigma(T)$ be an isolated eigenvalue with finite multiplicity $m$, and $M \subset \mathbb{C}$ an open neighborhood of $\lambda$ such that $\sigma(T) \cap M = \{\lambda\}$. Then:*

1. Convergence of eigenvalues: *There exists an $N \in \mathbb{N}$ such that, for all $n > N$, the set $\sigma(T_n) \cap M$ is an isolated part of $\sigma(T_n)$ consists of at most $m$ different eigenvalues, and their multiplicities sum up to $m$. Moreover, the sequence of sets $\sigma(T_n) \cap M$ converges to the set $\{\lambda\}$ in the sense that every sequence $(\lambda_n)_{n \in \mathbb{N}}$ with $\lambda_n \in \sigma(T_n) \cap M$ satisfies $\lim \lambda_n = \lambda$.*

2. Convergence of spectral projections: *Let $\mathrm{Pr}$ be the spectral projection of $T$ corresponding to $\lambda$, and for $n > N$, let $\mathrm{Pr}_n$ be the spectral projection of $T_n$ corresponding to $\sigma(T_n) \cap M$ (which is well defined according to part 1). Then $\mathrm{Pr}_n \xrightarrow{p} \mathrm{Pr}$.*

3. Convergence of eigenvectors: *If, additionally, $\lambda$ is a simple eigenvalue, then there exists some $N \in \mathbb{N}$ such that, for all $n > N$, the sets $\sigma(T_n) \cap M$ consist of a simple eigenvalue $\lambda_n$. The corresponding eigenfunctions $f_n$ converge up to a change of sign [i.e., there exists a sequence $(a_n)_n$ of signs $a_n \in \{-1, +1\}$ such that $a_n f_n$ converges].*

Proof. See Proposition 3.18 and Sections 3.6 and 5.1 in Chatelin [13]. $\square$

To prove rates of convergence, we will also need some quantitative perturbation theory results for spectral projections. The following theorem can be found in Atkinson [5]:

Theorem 7 (Atkinson [5]). *Let $(E, \|\cdot\|_E)$ be an arbitrary Banach space and $B$ its unit ball. Let $(K_n)_{n \in \mathbb{N}}$ and $K$ be compact linear operators on $E$ such that $K_n \xrightarrow{cc} K$. For a nonzero eigenvalue $\lambda \in \sigma(K)$, denote the corresponding spectral projection by $\mathrm{Pr}$. Let $M \subset \mathbb{C}$ be an open neighborhood of $\lambda$ such that $\sigma(K) \cap M = \{\lambda\}$. There exists some $N \in \mathbb{N}$ such that, for all $n > N$, the set $\sigma(K_n) \cap M$ is isolated in $\sigma(K_n)$. Let $\mathrm{Pr}_n$, the corresponding spectral projections of $K_n$ for $\sigma(K_n) \cap M$. Then there exists a constant $C > 0$ such that, for every $x \in \mathrm{Pr}\, E$,*

$$\|x - \mathrm{Pr}_n x\|_E \leq C(\|(K_n - K)x\|_E + \|x\|_E \|(K - K_n)K_n\|).$$

*The constant $C$ is independent of $x$, but it depends on $\lambda$ and $\sigma(K)$.*

For a probability measure $P$ and a function $f \in C(\mathcal{X})$, we introduce the abbreviation $Pf := \int f(x)\, dP(x)$. Let $(X_n)_n$ be a sequence of i.i.d. random variables drawn according to $P$, and denote by $P_n := 1/n \sum_{i=1}^n \delta_{X_i}$ the corresponding empirical distributions. A set $\mathcal{F} \subset C(\mathcal{X})$ is called a Glivenko–Cantelli class if

$$\sup_{f \in \mathcal{F}} |Pf - P_n f| \to 0 \qquad \text{a.s.}$$

Finally, the covering numbers $N(\mathcal{F}, \varepsilon, d)$ of a totally bounded set $\mathcal{F}$ with metric $d$ are defined as the smallest number $n$ such that $\mathcal{F}$ can be covered with $n$ balls of radius $\varepsilon$.

**5. Convergence of normalized spectral clustering.** In this section we present our results on the convergence of normalized spectral clustering. We start with an overview over our method, then prove several propositions, and finally state and prove our main theorems at the end of this section. The case of unnormalized spectral clustering will be treated in Section 7.

5.1. *Overview over the methods.* On a high level, the approach to prove convergence of spectral clustering is very similar in both the normalized and unnormalized case. In this section we focus on the normalized case. Moreover, as we have already seen that there is an explicit one-to-one relationship between the eigenvalues and eigenvectors of $L'_n$, $L''_n$ and the generalized eigenproblem $L_n v = \lambda D_n v$, we only consider the matrix $L'_n$ in the following. All results naturally can be carried over to the other cases. To study the convergence of spectral clustering, we have to investigate whether the eigenvectors of the Laplacians constructed on $n$ sample points "converge" for $n \to \infty$. For simplicity, let us discuss the case of the second eigenvector. For all $n \in \mathbb{N}$, let $v_n \in \mathbb{R}^n$ be the second eigenvector of $L'_n$. The technical difficulty for proving convergence of $(v_n)_{n \in \mathbb{N}}$ is that, for different sample sizes $n$, the vectors $v_n$ live in different spaces (as they have length $n$). Thus, standard notions of convergence cannot be applied. What we want to show instead is that there exists a function $f \in C(\mathcal{X})$ such that the difference between the eigenvector $v_n$ and the restriction of $f$ to the sample converges to 0, that is, $\|v_n - \rho_n f\|_\infty \to 0$. Our approach to achieve this takes one more detour. We replace the vector $v_n$ by a function $f_n \in C(\mathcal{X})$ such that $v_n = \rho_n f_n$. This function $f_n$ will be the second eigenfunction of an operator $U'_n$ acting on the space $C(\mathcal{X})$. Then we use the fact that

$$\|v_n - \rho_n f\|_\infty = \|\rho_n f_n - \rho_n f\|_\infty \leq \|f_n - f\|_\infty.$$

Hence, it will be enough to show that $\|f_n - f\|_\infty \to 0$. As the sequence, $f_n$ will be random, this convergence will hold almost surely.

STEP 1 [*Relating the matrices $L'_n$ to linear operators $U'_n$ on $C(\mathcal{X})$*]. First we will construct a family $(U'_n)_{n \in \mathbb{N}}$ of linear operators on $C(\mathcal{X})$ which, if restricted to the sample, "behaves" as $(L'_n)_{n \in \mathbb{N}}$: for all $f \in C(\mathcal{X})$, we will have the relation $\rho_n U'_n f = L'_n \rho_n f$. In the following we will then study the convergence of $(U'_n)_n$ on $C(\mathcal{X})$ instead of the convergence of $(L'_n)_n$.

STEP 2 [*Relation between $\sigma(L'_n)$ and $\sigma(U'_n)$*]. In Step 1 we replaced the *operators* $L'_n$ by operators $U'_n$ on $C(\mathcal{X})$. But as we are interested in the *eigenvectors* of $L'_n$, we have to check whether they can actually be recovered from the eigenfunctions of $U'_n$. By elementary linear algebra, we can prove that the "interesting" eigenfunctions $f_n$ and eigenvectors $v_n$ of $U'_n$ and $L'_n$ are in a one-to-one relationship and can be computed from each other by the relation $v_n = \rho_n f_n$. As a consequence, if the eigenfunctions $f_n$ of $U'_n$ converge, the same is true for the eigenvectors of $L'_n$.

STEP 3 (*Convergence of $U'_n \to U'$*). In this step we want to prove that certain eigenvalues and eigenfunctions of $U'_n$ converge to the corresponding quantities of some limit operator $U'$. For this, we will have to establish a rather strong type of convergence of linear operators. Pointwise convergence is in general too weak for this purpose; on the other hand, it will turn out that operator norm convergence does not hold in our context. The type of convergence we will consider is compact convergence, which is between pointwise convergence and operator norm convergence and is just strong enough for proving convergence of spectral properties. The notion of compact convergence has originally been developed in the context of (deterministic) numerical approximation of integral operators. We adapt those methods to a framework where the spectrum of a linear operator $U'$ is approximated by the spectra of *random* operators $U'_n$. Here, a key element is the fact that certain classes of functions are Glivenko–Cantelli classes: the integrals over the functions in those classes can be approximated uniformly by empirical integrals based on the random sample.

5.2. *Step* 1: *Construction of the operators on $C(\mathcal{X})$*. We define the following functions and operators, which are all supposed to act on $C(\mathcal{X})$: The degree functions

$$d_n(x) := \int k(x,y)\,dP_n(y) \in C(\mathcal{X}),$$

$$d(x) := \int k(x,y)\,dP(y) \in C(\mathcal{X}),$$

the multiplication operators,

$$M_{d_n} : C(\mathcal{X}) \to C(\mathcal{X}), \qquad M_{d_n}f(x) := d_n(x)f(x),$$

$$M_d : C(\mathcal{X}) \to C(\mathcal{X}), \qquad M_d f(x) := d(x)f(x),$$

the integral operators

$$S_n : C(\mathcal{X}) \to C(\mathcal{X}), \qquad S_n f(x) := \int k(x,y)f(y)\,dP_n(y),$$

$$S : C(\mathcal{X}) \to C(\mathcal{X}), \qquad S f(x) := \int k(x,y)f(y)\,dP(y),$$

and the corresponding differences

$$U_n : C(\mathcal{X}) \to C(\mathcal{X}), \qquad U_n f(x) := M_{d_n}f(x) - S_n f(x),$$

$$U : C(\mathcal{X}) \to C(\mathcal{X}), \qquad U f(x) := M_d f(x) - S f(x).$$

The operators $U_n$ and $U$ will be used to deal with the case of unnormalized spectral clustering. For the normalized case, we introduce the normalized similarity functions

$$h_n(x,y) := k(x,y)/\sqrt{d_n(x)d_n(y)},$$

$$h(x,y) := k(x,y)/\sqrt{d(x)d(y)},$$

the integral operators

$$T_n : C(\mathcal{X}) \to C(\mathcal{X}), \qquad T_n f(x) = \int h(x,y)f(y)\,dP_n(y),$$

$$\widehat{T_n} : C(\mathcal{X}) \to C(\mathcal{X}), \qquad \widehat{T_n} f(x) = \int h_n(x,y)f(y)\,dP_n(y),$$

$$T : C(\mathcal{X}) \to C(\mathcal{X}), \qquad T f(x) = \int h(x,y)f(y)\,dP(y),$$

and the differences

$$U'_n := I - \widehat{T_n},$$

$$U' := I - T.$$

In all what follows, the operators introduced above are always meant to act on the Banach space $(C(\mathcal{X}), \|\cdot\|_\infty)$, and their operator norms will be taken with respect to this space. We now summarize the properties of those operators in the following proposition. Recall the general assumptions and the definition of the restriction operator $\rho_n$ of Section 4.

PROPOSITION 8 (Relations between the operators).    *Under the general assumptions, the functions $d_n$ and $d$ are continuous, bounded from below by the constant $l > 0$, and from above by $\|k\|_\infty$. All operators defined above are bounded, and the integral operators are compact. The operator norms of $M_{d_n}$, $M_d$, $S_n$ and $S$ are bounded by $\|k\|_\infty$, the ones of $\widehat{T_n}$, $T_n$ and $T$ by $\|k\|_\infty/l$. Moreover, we have the following:*

$$\frac{1}{n} D_n \circ \rho_n = \rho_n \circ M_{d_n}, \qquad \frac{1}{n} K_n \circ \rho_n = \rho_n \circ S_n,$$

$$\frac{1}{n} L_n \circ \rho_n = \rho_n \circ U_n, \qquad L'_n \circ \rho_n = \rho_n \circ U'_n.$$

PROOF.    All statements follow directly from the definitions and the general assumptions. Note that in the case of the unnormalized Laplacian $L_n$ we get the scaling factor $1/n$ from the $1/n$-factor hidden in the empirical distribution $P_n$. In the case of the normalized Laplacian, this scaling factor cancels with the scaling factors of the degree functions in the denominators. □

The main statement of this proposition is that if restricted to the sample points, $U_n$ "behaves as" $\frac{1}{n}L_n$ and $U_n'$ as $L_n'$. Moreover, by the law of large numbers, it is clear that for fixed $f \in C(\mathcal{X})$ and $x \in \mathcal{X}$ the empirical quantities converge to the corresponding true quantities, in particular, $U_n f(x) \to U f(x)$ and $U_n' f(x) \to U' f(x)$. Proving stronger convergence statements will be the main part of Step 3.

5.3. *Step* 2: *Relations between the spectra.* The following proposition establishes the connections between the spectra of $L_n'$ and $U_n'$. We show that $U_n'$ and $L_n'$ have more or less the same spectrum and that the eigenfunctions $f$ of $U_n'$ and eigenvectors $v$ of $L_n'$ correspond to each other by the relation $v = \rho_n f$.

PROPOSITION 9 (Spectrum of $U_n'$).

1. *If $f \in C(\mathcal{X})$ is an eigenfunction of $U_n'$ with the eigenvalue $\lambda$, then the vector $v = \rho_n f \in \mathbb{R}^n$ is an eigenvector of the matrix $L_n'$ with eigenvalue $\lambda$.*

2. *Let $\lambda \neq 1$ be an eigenvalue of $U_n'$ with eigenfunction $f \in C(\mathcal{X})$, and $v := (v_1, \ldots, v_n) := \rho_n f \in \mathbb{R}^n$. Then $f$ is of the form*

$$(1) \qquad f(x) = \frac{1/n \sum_j k(x, X_j) v_j}{1 - \lambda}.$$

3. *If $v$ is an eigenvector of the matrix $L_n'$ with eigenvalue $\lambda \neq 1$, then $f$ defined by equation (1) is an eigenfunction of $U_n'$ with eigenvalue $\lambda$.*

4. *The spectrum of $U_n'$ consists of finitely many nonnegative eigenvalues with finite multiplicity. The essential spectrum of $U_n'$ consists of at most one point, namely, $\sigma_{\mathrm{ess}}(U_n') = \{1\}$. The spectrum of $U'$ consists of at most countably many nonnegative eigenvalues with finite multiplicity. Its essential spectrum consists at most of the point $\{1\}$, which is also the only possible accumulation point in $\sigma(U')$.*

PROOF. *Part* 1: It is obvious from Proposition 8 that $U_n' f = \lambda f$ implies $L_n' v = \lambda v$. Note also that part 2 shows that $v$ is not the constant 0 vector.

*Part* 2: Follows directly from solving the eigenvalue equation.

*Part* 3: Define $f$ as in equation (1). It is well defined because $v$ is an eigenvector of $\frac{1}{n}L_n$, and $f$ is an eigenfunction of $U_n$ with eigenvalue $\lambda$.

*Part* 4: According to Proposition 8, $\widehat{T_n}$ is a compact integral operator, and its essential spectrum is at most $\{0\}$. The spectrum $\sigma(U_n')$ of $U_n' = I - \widehat{T_n}$ is given by $1 - \sigma(\widehat{T_n})$. The statements about the eigenvalues of $U_n'$ follow from the properties of the eigenvalues of $L_n'$ and parts 1–3 of the proposition. An analogous reasoning leads to the statements for $U'$. □

This proposition establishes a one-to-one correspondence between the eigenvalues and eigenvectors of $L'_n$ and $U'_n$, provided they satisfy $\lambda \neq 1$. The condition $\lambda \neq 1$ needed to ensure that the denominator of equation (1) does not vanish. As a side remark, note that the set $\{1\}$ is the essential spectrum of $U'_n$. Thus, the condition $\lambda \neq 1$ can also be written as $\lambda \notin \sigma_{\mathrm{ess}}(U'_n)$, which will be analogous to the condition on the eigenvalues in the unnormalized case. This condition ensures that $\lambda$ is isolated in the spectrum.

5.4. *Step* 3: *Compact convergence.* In this section we want to prove that the sequence of random operators $U'_n$ converges compactly to $U'$ almost surely. First we will prove pointwise convergence. Note that on the space $C(\mathcal{X})$, the pointwise convergence of a sequence $U'_n$ of operators is defined as $\|U'_n f - U' f\|_\infty \to 0$, that is, for each $f \in C(\mathcal{X})$, the sequence $(U'_n f)_n$ has to converge uniformly over $\mathcal{X}$. To establish this convergence, we will need to show that several classes of functions are "not too large," that is, they are Glivenko–Cantelli classes. For convenience, we introduce the following notation:

DEFINITION 10 (*Particular sets of functions*).    Let $k: \mathcal{X} \times \mathcal{X} \to \mathbb{R}$ be a similarity function, $h: \mathcal{X} \times \mathcal{X} \to \mathbb{R}$ the corresponding normalized similarity function as introduced above and $g \in C(\mathcal{X})$ an arbitrary function. We use the shorthand notation $k(x, \cdot)$, $g(\cdot)k(x, \cdot)$ and $h(x, \cdot)h(y, \cdot)$ to denote the functions $z \mapsto k(x, z)$, $z \mapsto g(z)k(x, z)$ and $z \mapsto h(x, z)h(y, z)$. We define the following:

$$\mathcal{K} := \{k(x, \cdot); x \in \mathcal{X}\}, \qquad \mathcal{H} := \{h(x, \cdot); x \in \mathcal{X}\},$$
$$g \cdot \mathcal{H} := \{g(\cdot)h(x, \cdot); x \in \mathcal{X}\}, \qquad \mathcal{H} \cdot \mathcal{H} := \{h(x, \cdot)h(y, \cdot); x, y \in \mathcal{X}\}.$$

PROPOSITION 11 (Glivenko–Cantelli classes).    *Under the general assumptions, the classes $\mathcal{K}$, $\mathcal{H}$ and $g \cdot \mathcal{H}$ [for arbitrary $g \in C(\mathcal{X})$] are Glivenko–Cantelli classes.*

PROOF.    As $k$ is a continuous function defined on a compact domain, it is uniformly continuous. In this case it is easy to construct, for each $\varepsilon > 0$, a finite $\varepsilon$-cover with respect to $\|\cdot\|_\infty$ of $\mathcal{K}$ from a finite $\delta$-cover of $\mathcal{X}$. Hence, $\mathcal{K}$ has finite $\|\cdot\|_\infty$-covering numbers. Then it is easy to see that $\mathcal{K}$ also has finite $\|\cdot\|_{L_1(P)}$-bracketing numbers (cf. van der Vaart and Wellner [45], page 84). Now the statement about the class $\mathcal{K}$ follows from Theorem 2.4.1 of van der Vaart and Wellner [45]. The statements about the classes $\mathcal{H}$ and $g \cdot \mathcal{H}$ can be proved in the same way, hereby observing that $h$ is continuous and bounded as a consequence of the general assumptions.    □

Note that it is a direct consequence of this proposition that the empirical degree function $d_n$ converges uniformly to the true degree function $d$, that is,

$$\|d_n - d\|_\infty = \sup_{x \in \mathcal{X}} |d_n(x) - d(x)| = \sup_{x \in \mathcal{X}} |P_n k(x, \cdot) - P k(x, \cdot)| \to 0 \qquad \text{a.s.}$$

PROPOSITION 12 ($\widehat{T_n}$ converges pointwise to $T$ a.s.).   *Under the general assumptions, $\widehat{T_n} \xrightarrow{p} T$ almost surely.*

PROOF.   For arbitrary $f \in C(\mathcal{X})$, we have

$$\|\widehat{T_n} f - Tf\|_\infty \leq \|\widehat{T_n} f - T_n f\|_\infty + \|T_n f - Tf\|_\infty.$$

The second term can be written as

$$\|T_n f - Tf\|_\infty = \sup_{x \in \mathcal{X}} |P_n(h(x, \cdot) f(\cdot)) - P(h(x, \cdot) f(\cdot))| = \sup_{g \in f \cdot \mathcal{H}} |P_n g - P g|,$$

which converges to 0 a.s. by Proposition 11. The first term can be bounded by

$$\|T_n f - \widehat{T_n} f\|_\infty \leq \|f\|_\infty \|k\|_\infty \sup_{x,y \in \mathcal{X}} \left| \frac{1}{\sqrt{d_n(x) d_n(y)}} - \frac{1}{\sqrt{d(x) d(y)}} \right|$$

$$= \|f\|_\infty \frac{\|k\|_\infty}{l^2} \sup_{x,y \in \mathcal{X}} \frac{|d_n(x) d_n(y) - d(x) d(y)|}{\sqrt{d_n(x) d_n(y)} + \sqrt{d(x) d(y)}}$$

$$\leq \|f\|_\infty \frac{\|k\|_\infty}{2l^3} \sup_{x,y \in \mathcal{X}} |d_n(x) d_n(y) - d(x) d(y)|$$

$$\leq \|f\|_\infty \frac{\|k\|_\infty^2}{l^3} |d_n(x) - d(x)| \leq \|f\|_\infty \frac{\|k\|_\infty^2}{l^3} \sup_{g \in \mathcal{K}} |P_n g - P g|.$$

Together with Proposition 11 this finishes the proof.   □

PROPOSITION 13 ($\widehat{T_n}$ converges collectively compactly to $T$ a.s.).   *Under the general assumptions, $\widehat{T_n} \xrightarrow{cc} T$ almost surely.*

PROOF.   We have already seen the pointwise convergence $\widehat{T_n} \xrightarrow{p} T$ in Proposition 12. Next we have to prove that, for some $N \in \mathbb{N}$, the sequence of operators $(\widehat{T_n} - T)_{n > N}$ is collectively compact a.s. As $T$ is compact itself, it is enough to show that $(\widehat{T_n})_{n > N}$ is collectively compact a.s. This will be done using the Arzela–Ascoli theorem (e.g., Section I.6 of Reed and Simon [41]). First we fix the random sequence $(X_n)_n$ and, hence, the random operators $(\widehat{T_n})_n$. By Proposition 8, we know that $\|\widehat{T_n}\| \leq \|k\|_\infty / l$ for all $n \in \mathbb{N}$. Hence, the functions in $\bigcup_n \widehat{T_n}B$ are uniformly bounded by $\sup_{n \in \mathbb{N}, f \in B} \|\widehat{T_n}f\|_\infty \leq \|k\|_\infty/l$. To prove that the functions in $\bigcup_{n>N} \widehat{T_n}B$ are equicontinuous, we have to bound the expression $|g(x) - g(x')|$ in terms of the distance between $x$ and $x'$, uniformly in $g \in \bigcup_n \widehat{T_n}B$. For fixed sequence $(X_n)_{n \in \mathbb{N}}$ and all $n \in \mathbb{N}$, we have that for all $x, x' \in \mathcal{X}$,

$$
\sup_{f \in B, n \in \mathbb{N}} |\widehat{T_n}f(x) - \widehat{T_n}f(x')| = \sup_{f \in B, n \in \mathbb{N}} \left| \int (h_n(x,y) - h_n(x',y))f(y)\, dP_n(y) \right|
$$
$$
\leq \sup_{f \in B, n \in \mathbb{N}} \|f\|_\infty \int |h_n(x,y) - h_n(x',y)|\, dP_n(y)
$$
$$
\leq \|h_n(x, \cdot) - h_n(x', \cdot)\|_\infty.
$$

Now we have to prove that the right-hand side gets small whenever the distance between $x$ and $x'$ gets small:

$$
\sup_y |h_n(x,y) - h_n(x',y)|
$$
$$
\leq \frac{1}{l^{3/2}}(\|\sqrt{d_n}\|_\infty \|k(x,\cdot) - k(x',\cdot)\|_\infty + \|k\|_\infty |\sqrt{d_n(x)} - \sqrt{d_n(x')}|)
$$
$$
\leq \frac{1}{l^{3/2}}\left( \|k\|_\infty^{1/2}\|k(x,\cdot) - k(x',\cdot)\|_\infty + \frac{\|k\|_\infty}{2l^{1/2}}|d_n(x) - d_n(x')| \right)
$$
$$
\leq C_1\|k(x,\cdot) - k(x',\cdot)\|_\infty + C_2|d(x) - d(x')| + C_3\|d_n - d\|_\infty.
$$

As $\mathcal{X}$ is a compact space, the continuous functions $k$ (on the compact space $\mathcal{X} \times \mathcal{X}$) and $d$ are in fact uniformly continuous. Thus, the first two (deterministic) terms $\|k(x,\cdot) - k(x',\cdot)\|_\infty$ and $|d(x) - d(x')|$ can be made arbitrarily small for all $x, x'$ whenever the distance between $x$ and $x'$ is small. For the third term $\|d_n - d\|_\infty$, which is a random term, we know by the Glivenko–Cantelli properties of Proposition 11 that it converges to 0 a.s. This means that for each given $\varepsilon > 0$ there exists some $N \in \mathbb{N}$ such that, for all $n > N$, we have $\|d_n - d\|_\infty \leq \varepsilon$ a.s. Together, these arguments show that $\bigcup_{n>N} \widehat{T_n}B$ is equicontinuous a.s. By the Arzela–Ascoli theorem, we then know that $\bigcup_{n>N} \widehat{T_n}B$ is relatively compact a.s., which concludes the proof. $\square$

PROPOSITION 14 ($U_n'$ converges compactly to $U'$ a.s.).    *Under the general assumptions, $U_n' \xrightarrow{c} U'$ a.s.*

PROOF.    This follows directly from the facts that collectively compact convergence implies compact convergence, the definitions of $U_n'$ to $U'$, and Proposition 13. $\square$

5.5. *Assembling all pieces.* Now we have collected all ingredients to state and prove our convergence result for normalized spectral clustering. The following theorem is the precisely formulated version of the informal Result 1 of the introduction:

THEOREM 15 (Convergence of normalized spectral clustering). *Assume that the general assumptions hold. Let $\lambda \neq 1$ be an eigenvalue of $U'$ and $M \subset \mathbb{C}$ an open neighborhood of $\lambda$ such that $\sigma(U') \cap M = \{\lambda\}$. Then:*

1. Convergence of eigenvalues: *The eigenvalues in $\sigma(L'_n) \cap M$ converge to $\lambda$ in the sense that every sequence $(\lambda_n)_{n \in \mathbb{N}}$ with $\lambda_n \in \sigma(L'_n) \cap M$ satisfies $\lambda_n \to \lambda$ almost surely.*
2. Convergence of spectral projections: *There exists some $N \in \mathbb{N}$ such that, for $n > N$, the sets $\sigma(U'_n) \cap M$ are isolated in $\sigma(U'_n)$. For $n > N$, let $\mathrm{Pr}'_n$ be the spectral projections of $U'_n$ corresponding to $\sigma(U'_n) \cap M$, and $\mathrm{Pr}$ the spectral projection of $U$ for $\lambda$. Then $\mathrm{Pr}'_n \xrightarrow{p} \mathrm{Pr}$ a.s.*
3. Convergence of eigenvectors: *If $\lambda$ is a simple eigenvalue, then the eigenvectors of $L'_n$ converge a.s. up to a change of sign: if $v_n$ is the eigenvector of $L'_n$ with eigenvalue $\lambda_n$, $v_{n,i}$ its ith coordinate, and $f$ the eigenfunction of eigenvalue $\lambda$, then there exists a sequence $(a_n)_{n \in \mathbb{N}}$ with $a_i \in \{+1, -1\}$ such that $\sup_{i=1,\dots,n} |a_n v_{n,i} - f(X_i)| \to 0$ a.s. In particular, for all $b \in \mathbb{R}$, the sets $\{a_n f_n > b\}$ and $\{f > b\}$ converge, that is, their symmetric difference satisfies $P(\{f > b\} \triangle \{a_n f_n > b\}) \to 0$.*

PROOF. In Proposition 9 we established a one-to-one correspondence between the eigenvalues $\lambda \neq 1$ of $L'_n$ and $U'_n$, and we saw that the eigenvalues $\lambda$ of $U'$ with $\lambda \neq 1$ are isolated and have finite multiplicity. In Proposition 14 we proved the compact convergence of $U'_n$ to $U'$, which according to Proposition 6 implies the convergence of the spectral projections of isolated eigenvalues with finite multiplicity. For simple eigenvalues, this implies the convergence of the eigenvectors up to a change of sign. The convergence of the sets $\{f_n > b\}$ is a simple consequence of the almost sure convergence of $(a_n f_n)_n$. □

Observe that we only get convergence of the *eigenvectors* if the eigenvalue of the limit operator is simple. If this assumption is not satisfied, we only get convergence of the *eigenspaces*, but not of the individual eigenvectors.

**6. Rates of convergence in the normalized case.** In this section we want to prove statements about the rates of convergence of normalized spectral clustering. Our main result is the following:

THEOREM 16 (Rate of convergence of normalized spectral clustering). *Under the general assumptions, let $\lambda \neq 0$ be a simple eigenvalue of $T$ with eigenfunction $u$, $(\lambda_n)_n$ a sequence of eigenvalues of $\widehat{T_n}$ such that $\lambda_n \to \lambda$, and $(u_n)_n$ a corresponding sequence of eigenfunctions. Define $\mathcal{F} = \mathcal{K} \cup (u \cdot \mathcal{H}) \cup (\mathcal{H} \cdot \mathcal{H})$. Then there exists a constant $C' > 0$ [which only depends on the similarity function $k$, on $\sigma(T)$ and on $\lambda$] and a sequence $(a_n)_n$ of signs $a_n \in \{+1, -1\}$ such that*

$$\|a_n u_n - u\|_\infty \leq C' \sup_{f \in \mathcal{F}} |P_n f - P f|.$$

This theorem shows that the rate of convergence of normalized spectral clustering is at least as good as the rate of convergence of the supremum of the empirical process indexed by $\mathcal{F}$. To determine the latter, there exist a variety of tools and techniques from the theory of empirical processes, such as covering numbers, VC dimension and Rademacher complexities; see, for example, van der Vaart and Wellner [45], Dudley [18], Mendelson [34] and Pollard [39]. In particular, it is the case that "the nicer" the kernel function $k$ is (e.g., $k$ is Lipschitz, or smooth, or positive definite), the faster the rate of convergence on the right-hand side will be. As an example we will consider the case of the Gaussian similarity function $k(x, y) = \exp(-\|x - y\|^2 / \sigma^2)$, which is widely used in practical applications of spectral clustering.

EXAMPLE 1 (Rate of convergence for Gaussian kernel). *Let $\mathcal{X}$ be compact subset of $\mathbb{R}^d$ and $k(x, y) = \exp(-\|x - y\|^2 / \sigma^2)$. Then the eigenvectors in Theorem 16 converge with rate $\mathcal{O}(1/\sqrt{n})$.*

For the case of unnormalized spectral clustering, it is possible to obtain similar results on the speed of convergence, for example, by using Proposition 5.3 in Chapter 5 of Chatelin [13] instead of the results of Atkinson [5] (note that in the unnormalized case, the assumptions of Theorem 7 are not satisfied, as we only have compact convergence instead of collectively compact convergence). As we recommend to use normalized rather than unnormalized spectral clustering anyway, we do not discuss this issue any further. The remaining part of this section is devoted to the proofs of Theorem 16 and Example 1.

6.1. *Some technical preparations.* Before we can prove Theorem 16 we need to show several technical propositions.

PROPOSITION 17 (Some technical bounds). *Assume that the general conditions are satisfied, and let $g \in CX$. Then the following bounds hold*

*true:*

$$\|\widehat{T_n} - T_n\| \leq \frac{\|k\|_\infty^2}{l^3} \sup_{f \in \mathcal{K}} |P_n f - P f|,$$

$$\|(T_n - T)g\|_\infty \leq \sup_{f \in g \cdot \mathcal{H}} |P_n f - P f|,$$

$$\|(T - T_n)T_n\| \leq \sup_{f \in \mathcal{H} \cdot \mathcal{H}} |P_n f - P f|.$$

PROOF. The first inequality can be proved similarly to Proposition 12, the second inequality is a direct consequence of the definitions. The third inequality follows by straight forward calculations similar to the ones in the previous section and using Fubini's theorem and the symmetry of $h$. $\quad\square$

PROPOSITION 18 (Convergence of one-dimensional projections). *Let* $(v_n)_n$ *be a sequence of vectors in some Banach space* $(E, \|\cdot\|)$ *with* $\|v_n\| = 1$, $\mathrm{Pr}_n$ *the projections on the one-dimensional subspace spanned by* $v_n$, *and* $v \in E$ *with* $\|v\| = 1$. *Then there exists a sequence* $(a_n)_n \in \{+1, -1\}$ *of signs such that*

$$\|a_n v_n - v\| \leq 2\|v - \mathrm{Pr}_n v\|.$$

*In particular, if* $\|v - \mathrm{Pr}_n v\| \to 0$, *then* $v_n$ *converges to* $v$ *up to a change of sign.*

PROOF. By the definition of $\mathrm{Pr}_n$, we know that $\mathrm{Pr}_n v = c_n v_n$ for some $c_n \in \mathbb{R}$. Define $a_n := \mathrm{sgn}(c_n)$. Then

$$|a_n - c_n| = |1 - |c_n|| = |\|v\| - |c_n| \cdot \|v_n\|| \leq \|v - c_n v_n\| = \|v - \mathrm{Pr}_n v\|.$$

From this, we can conclude that

$$\|v - a_n v_n\| \leq \|v - c_n v_n\| + \|c_n v_n - a_n v_n\| \leq 2\|v - \mathrm{Pr}_n v\|. \quad\square$$

6.2. *Proof of Theorem 16.* First we fix a realization of the random variables $(X_n)_n$. From the convergence of the spectral projections in Theorem 15 we know that if $\lambda \in \sigma(T)$ is simple, so are $\lambda_n \in \sigma(\widehat{T_n})$ for large $n$. Then the eigenfunctions $u_n$ are uniquely determined up to a change of orientation. In Proposition 18 we have seen that the speed of convergence of $u_n$ to $u$ is bounded by the speed of convergence of the expression $\|u - \mathrm{Pr}_n u\|$ from Theorem 7. As we already know by Section 5, the operators $\widehat{T_n}$ and $T$ satisfy the assumptions in Theorem 7. Accordingly, $\|u - \mathrm{Pr}_n u\|$ can be bounded by the two terms $\|(\widehat{T_n} - T)u\|$ and $\|(T - \widehat{T_n})\widehat{T_n}\|$. It will turn out that both terms are easier to bound if we can replace the operator $\widehat{T_n}$ by $T_n$. To accomplish this, observe that

$$\|(T - \widehat{T_n})\widehat{T_n}\| \le \|T\|\|T_n - \widehat{T_n}\| + \|(T - T_n)T_n\|$$
$$+ \|T_nT_n - T_n\widehat{T_n}\| + \|T_n\widehat{T_n} - \widehat{T_n}\widehat{T_n}\|$$
$$\le 3\frac{\|k\|_\infty}{l}\|T_n - \widehat{T_n}\| + \|(T - T_n)T_n\|$$

and also

$$\|(\widehat{T_n} - T)u\|_\infty \le \|u\|_\infty\|\widehat{T_n} - T_n\| + \|(T_n - T)u\|_\infty.$$

Note that $T_n$ does not converge to $T$ in operator norm (cf. page 197 in Section 4.7.4 of Chatelin [13]). Thus, it does not make sense to bound $\|(T_n - T)u\|_\infty$ by $\|T_n - T\|\|u\|_\infty$ or $\|(T - T_n)T_n\|$ by $\|T - T_n\|\|T_n\|$. Assembling all inequalities, applying Proposition 18 and Theorem 7, and choosing the signs $a_n$ as in the proof of Proposition 18, we obtain

$$\|a_nu_n - u\| \le 2\|u - \mathrm{Pr}_{\lambda_n}u\| \le 2C(\|(\widehat{T_n} - T)u\| + \|(T - \widehat{T_n})\widehat{T_n}\|)$$
$$\le 2C\left(\left(\frac{3\|k\|_\infty}{l} + 1\right)\|T_n - \widehat{T_n}\| + \|(T_n - T)u\|_\infty + \|(T - T_n)T_n\|\right)$$
$$\le C' \sup_{f\in\mathcal{K}\cup(u\cdot\mathcal{H})\cup(\mathcal{H}\cdot\mathcal{H})} |P_nf - Pf|.$$

Here the last step was obtained by applying Proposition 17 and merging all occurring constants to one larger constant $C'$. As all arguments hold for each fixed realization $(X_n)_n$ of the sample points, they also hold for the random variables themselves almost surely. This concludes the proof of Theorem 16.

6.3. *Rate of convergence for the Gaussian kernel.* In this subsection we want to prove the convergence rate $\mathcal{O}(1/\sqrt{n})$ stated in Example 1 for the case of a Gaussian kernel function $k(x, y) = \exp(-\|x - y\|^2/\sigma^2)$. In principle, there are many ways to compute rates of convergence for terms of the form $\sup_f |Pf - P_nf|$ (see, e.g., van der Vaart and Wellner [45]). As discussing those methods is not the main focus of our paper, we choose a rather simple covering number approach which suffices for our purposes, but might not lead to the sharpest possible bounds. We will use the following theorem, which is well known in empirical process theory (nevertheless, we did not find a good reference for it; it can be obtained for example by combining Section 3.4 of Anthony [4], and Theorem 2.34 in Mendelson [34]):

THEOREM 19 (Entropy bound). *Let $(\mathcal{X}, \mathcal{A}, P)$ be an arbitrary probability space, $\mathcal{F}$ a class of real-valued functions on $\mathcal{X}$ with $\|f\|_\infty \le 1$. Let $(X_n)_{n\in\mathbb{N}}$ be a sequence of i.i.d. random variables drawn according to $P$, and $(P_n)_{n\in\mathbb{N}}$* *the corresponding empirical distributions. Then there exists some constant $c > 0$ such that, for all $n \in \mathbb{N}$ with probability at least $1 - \delta$,*

$$\sup_{f \in \mathcal{F}} |P_n f - P f| \le \frac{c}{\sqrt{n}} \int_0^\infty \sqrt{\log N(\mathcal{F}, \varepsilon, L_2(P_n))} \, d\varepsilon + \sqrt{\frac{1}{2n} \log \frac{2}{\delta}}.$$

We can see that if $\int_0^\infty \sqrt{\log N(\mathcal{F}, \varepsilon, L_2(P_n))} \, d\varepsilon < \infty$, then the whole expression scales as $\mathcal{O}(1/\sqrt{n})$. As a first step we would like to evaluate this integral for the function class $\mathcal{F} := \mathcal{K}$. To this, end we use bounds on the $\| \cdot \|_\infty$-covering numbers of $\mathcal{K}$ obtained in Proposition 1 in [51]. There it was proved that for $\varepsilon < c_0$ for a certain constant $c_0 > 0$ only depending to the kernel width $\sigma$, and for some constant $C$ which just depends on the dimension of the underlying space, the covering numbers satisfy

$$\log N(\mathcal{K}, \varepsilon, \| \cdot \|_\infty) \le C \left( \log \frac{1}{\varepsilon} \right)^2.$$

Plugging this into the integral, above we get

$$\int_0^\infty \sqrt{\log N(\mathcal{K}, \varepsilon, L_2(P_n))} \, d\varepsilon$$

$$\le \int_0^2 \sqrt{\log N(\mathcal{K}, \varepsilon, \| \cdot \|_\infty)} \, d\varepsilon$$

$$\le \sqrt{C} \int_0^{c_0} \log \frac{1}{\varepsilon} \, d\varepsilon + \int_{c_0}^2 \sqrt{\log N(\mathcal{K}, \varepsilon, \| \cdot \|_\infty)} \, d\varepsilon$$

$$\le \sqrt{C} c_0 (1 - \log c_0) + (2 - c_0) \sqrt{\log N(\mathcal{K}, c_0, \| \cdot \|_\infty)} < \infty.$$

According to Theorem 16, we have to use the entropy bound not only for the function class $\mathcal{F} = \mathcal{K}$, but for the class $\mathcal{F} = \mathcal{K} \cup (u \cdot \mathcal{H}) \cup (\mathcal{H} \cdot \mathcal{H})$. To this end, we will bound the $\| \cdot \|_\infty$-covering numbers of $\mathcal{K} \cup (u \cdot \mathcal{H}) \cup (\mathcal{H} \cdot \mathcal{H})$ in terms of the covering numbers of $\mathcal{K}$.

PROPOSITION 20 (Covering numbers). *Under the general assumptions, the following covering number bounds hold true:*

$$N(\mathcal{H}, \varepsilon, \| \cdot \|_\infty) \le N(\mathcal{K}, s\varepsilon, \| \cdot \|_\infty),$$

$$N(\mathcal{K} \cup (u \cdot \mathcal{H}) \cup (\mathcal{H} \cdot \mathcal{H}), \varepsilon, \| \cdot \|_\infty) \le 3 N(\mathcal{K}, q\varepsilon, \| \cdot \|_\infty),$$

*where* $s = \frac{\|k\|_\infty + 2\sqrt{l\|k\|_\infty}}{2l^2}$, $q := \min\{1, \|u\|_\infty s, \frac{\|k\|_\infty}{l} s\}$ *and* $u \in C(\mathcal{X})$ *arbitrary.*

This can be proved by straight forward calculations similar to the ones presented in the previous sections.

Combining this proposition with the integral bound for the Gaussian kernel as computed above, we obtain

$$\int_0^\infty \sqrt{\log N(\mathcal{F}, \varepsilon, L_2(P_n))}\, d\varepsilon \leq \int_0^\infty \sqrt{\log 3N(\mathcal{K}, q\varepsilon, \|\cdot\|_\infty)}\, d\varepsilon < \infty.$$

The entropy bound in Theorem 19 hence shows that the rate of convergence of $\sup_{f \in \mathcal{F}} |P_n f - Pf|$ is $\mathcal{O}(1/\sqrt{n})$, and by Theorem 16, the same now holds for the eigenfunctions of normalized spectral clustering.

**7. The unnormalized case.** Now we want to turn our attention to the case of unnormalized spectral clustering. It will turn out that this case is not as nice as the normalized case, as the convergence results will hold under strong conditions only. Moreover, those conditions are often violated in practice. In this case, the eigenvectors do not contain any useful information about the clustering of the data space.

7.1. *Convergence of unnormalized spectral clustering.* The main theorem about convergence of unnormalized spectral clustering (which was informally stated as Result 2 in Section 3) is as follows:

THEOREM 21 (Convergence of unnormalized spectral clustering). *Assume that the general assumptions hold. Let $\lambda \notin \mathrm{rg}(d)$ be an eigenvalue of $U$ and $M \subset \mathbb{C}$ an open neighborhood of $\lambda$ such that $\sigma(U) \cap M = \{\lambda\}$. Then:*

1. Convergence of eigenvalues: *The eigenvalues in $\sigma(\frac{1}{n} L_n) \cap M$ converge to $\lambda$ in the sense that every sequence $(\lambda_n)_{n \in \mathbb{N}}$ with $\lambda_n \in \sigma(\frac{1}{n} L_n) \cap M$ satisfies $\lambda_n \to \lambda$ almost surely.*

2. Convergence of spectral projections: *There exists some $N \in \mathbb{N}$ such that, for $n > N$, the sets $\sigma(U_n) \cap M$ are isolated in $\sigma(U_n)$. For $n > N$, let $\mathrm{Pr}_n$ be the spectral projections of $U_n$ corresponding to $\sigma(U_n) \cap M$, and $\mathrm{Pr}$ the spectral projection of $U$ for $\lambda$. Then $\mathrm{Pr}_n \xrightarrow{p} \mathrm{Pr}$ a.s.*

3. Convergence of eigenvectors: *If $\lambda$ is a simple eigenvalue, then the eigenvectors of $\frac{1}{n} L_n$ converge a.s. up to a change of sign: if $v_n$ is the eigenvector of $\frac{1}{n} L_n$ with eigenvalue $\lambda_n$, $v_{n,i}$ its $i$th coordinate, and $f$ the eigenfunction of $U$ with eigenvalue $\lambda$, then there exists a sequence $(a_n)_{n \in \mathbb{N}}$ with $a_i \in \{+1, -1\}$ such that $\sup_{i=1,\ldots,n} |a_n v_{n,i} - f(X_i)| \to 0$ a.s. In particular, for all $b \in \mathbb{R}$, the sets $\{a_n f_n > b\}$ and $\{f > b\}$ converge, that is, their symmetric difference satisfies $P(\{f > b\} \triangle \{a_n f_n > b\}) \to 0$.*

This theorem looks very similar to Theorem 15. The only difference is that the condition $\lambda \neq 1$ of Theorem 15 is now replaced by $\lambda \notin \mathrm{rg}(d)$. Note that in both cases, those conditions are equivalent to saying that $\lambda$ must be an isolated eigenvalue. In the normalized case, this is satisfied for all eigenvalues but $\lambda = 1$, as $U' = I - T'$ where $T'$ is a compact operator. In the unnormalized case, however, this condition can be violated, as the spectrum of $U$ contains a large continuous spectrum. Later we will see that this indeed leads to serious problems.

The proof of Theorem 7 is very similar to the one we presented in Section 5. The main difference between both cases is the structure of the spectra of $U_n$ and $U$. The proposition corresponding to Proposition 9 is the following:

PROPOSITION 22 (Spectrum of $U_n$).

1. If $f \in C(\mathcal{X})$ is an eigenfunction of $U_n$ with arbitrary eigenvalue $\lambda$, then the vector $v = \rho_n f \in \mathbb{R}^n$ is an eigenvector of the matrix $\frac{1}{n} L_n$ with eigenvalue $\lambda$.

2. Let $\lambda \notin \mathrm{rg}(d_n)$ be an eigenvalue of $U_n$ with eigenfunction $f \in C(\mathcal{X})$, and $v := (v_1, \ldots, v_n) := \rho_n f \in \mathbb{R}^n$. Then $f$ is of the form

$$(2) \qquad f(x) = \frac{1/n \sum_j k(x, X_j) v_j}{d_n(x) - \lambda}.$$

3. If $v$ is an eigenvector of the matrix $\frac{1}{n} L_n$ with eigenvalue $\lambda \notin \mathrm{rg}(d_n)$, then $f$ defined by equation (2) is an eigenfunction of $U_n$ with eigenvalue $\lambda$.

4. The essential spectrum of $U_n$ coincides with the range of the degree function, that is, $\sigma_{\mathrm{ess}}(U_n) = \mathrm{rg}(d_n)$. All eigenvalues of $U_n$ are nonnegative and can have accumulation points only in $\mathrm{rg}(d_n)$. The analogous statements also hold for the operator $U$.

PROOF. The first parts can be proved analogously to Proposition 9. For the last part, remember that the essential spectrum of the multiplication operator $M_{d_n}$ consists of the range of the multiplier function $d_n$. As $S_n$ is a compact operator, the essential spectrum of $U_n = M_{d_n} - S_n$ coincides with the essential spectrum of $M_{d_n}$, as we have already mentioned in the beginning of Section 4. The accumulation points of the spectrum of a bounded operator always belong to the essential spectrum. Finally, to see the nonnegativity of the eigenvalues, observe that if we consider the operator $U_n$ as an operator on $L_2(P_n)$ we have

$$\langle U_n f, f \rangle = \int \int (f(x) - f(y)) f(x) k(x, y) \, dP_n(y) \, dP_n(x)$$
$$= \frac{1}{2} \int \int (f(x) - f(y))^2 k(x, y) \, dP_n(y) \, dP_n(x) \geq 0.$$

Thus, $U$ is a nonnegative operator on $L_2(P_n)$ and as such only has a nonnegative eigenvalues. As we have $C(\mathcal{X}) \subset L_2(P)$ by the compactness of $\mathcal{X}$, the same holds for the eigenvalues of $U$ as an operator on $C(\mathcal{X})$. $\square$

This proposition establishes a one-to-one relationship between the eigenvalues of $U_n$ and $\frac{1}{n}L_n$, provided the condition $\lambda \notin \mathrm{rg}(d_n)$ is satisfied. Next we need to prove the compact convergence of $U_n$ to $U$:

PROPOSITION 23 ($U_n$ converges compactly to $U$ a.s.).    *Under the general assumptions, $U_n \xrightarrow{c} U$ a.s.*

PROOF. We consider the multiplication and integral operator parts of $U_n$ separately. Similarly to Proposition 13, we can prove that the integral operators $S_n$ converge collectively compactly to $S$ a.s., and, as a consequence, also $S_n \xrightarrow{c} S$ a.s. For the multiplication operators, we have operator norm convergence

$$\|M_{d_n} - M_d\| = \sup_{\|f\|_\infty \leq 1} \|d_n f - df\|_\infty \leq \|d_n - d\|_\infty \to 0 \qquad \text{a.s.}$$

by the Glivenko–Cantelli Proposition 11. As operator norm convergence implies compact convergence, we also have $M_{d_n} \xrightarrow{c} M_d$ a.s. Finally, it is easy to see that the sum of two compactly converging operators also converges compactly. □

Now Theorem 21 follows by a proof similar to the one of Theorem 15.

**8. Nonisolated eigenvalues.** The most important difference between the limit operators of normalized and unnormalized spectral clustering is the condition under which eigenvalues of the limit operator are isolated in the spectrum. In the normalized case this is true for all eigenvalues $\lambda \neq 1$, while in the unnormalized case this is only true for all eigenvalues satisfying $\lambda \notin \mathrm{rg}(d)$. In this section we want to investigate those conditions more closely. We will see that, especially in the unnormalized case, this condition can be violated, and that in this case spectral clustering will not yield sensible results. In particular, the condition $\lambda \notin \mathrm{rg}(d)$ is not an artifact of our methods, but plays a fundamental role. It is the main reason why we suggest to use normalized rather than unnormalized spectral clustering.

8.1. *Theoretical results.* First we will construct a simple example where all nontrivial eigenvalues $\lambda_2, \lambda_3, \ldots$ lie inside the range of the degree function.

EXAMPLE 2 [$\lambda_2 \notin \mathrm{rg}(d)$ violated].    *Consider the data space $\mathcal{X} = [1, 2] \subset \mathbb{R}$ and the probability distribution given by a piecewise constant probability density function $p$ on $\mathcal{X}$ with $p(x) = s$ if $4/3 \leq x < 5/3$ and $p(x) = (3 - s)/2$ otherwise, for some fixed constant $s \in [0, 3]$ (for example, for $s = 0.3$, this density has two clearly separated high density regions). As similarity function, we choose $k(x, y) := xy$. Then the only eigenvalue of $U$ outside of $\mathrm{rg}(d)$ is the trivial eigenvalue 0 with multiplicity one.*

PROOF. In this example, it is straightforward to verify that the degree function is given as $d(x) = 1.5x$ (independently of $s$) and has range $[1.5, 3]$ on $\mathcal{X}$. A function $f \in C(\mathcal{X})$ is an eigenfunction with eigenvalue $\lambda \notin \mathrm{rg}(d)$ of $U$ if the eigenvalue equation is satisfied:

$$(3) \qquad Uf(x) = d(x)f(x) - x \int yf(y)p(y)\,dy = \lambda f(x).$$

Defining the real number $\beta := \int yf(y)p(y)\,dy$, we can solve equation (3) for $f(x)$ to obtain $f(x) = \frac{\beta x}{d(x) - \lambda}$. Plugging this into the definition of $\beta$ yields the condition

$$(4) \qquad 1 = \int \frac{y^2}{d(y) - \lambda} p(y)\,dy.$$

Hence, $\lambda$ is an eigenvalue of $U$ if equation (4) is satisfied. For our simple density function $p$, the integral in this condition can be solved analytically. It can then be seen that $g(\lambda) := \int \frac{y^2}{d(y) - \lambda} p(y)\,dy = 1$ is only satisfied for $\lambda = 0$, hence, the only eigenvalue outside of $\mathrm{rg}(d)$ is the trivial eigenvalue 0 with multiplicity one. □

In the above example we can see that there indeed exist situations where the operator $U$ does not possess a nonzero eigenvalue with $\lambda \notin \mathrm{rg}(d)$. The next question is what happens in this situation.

PROPOSITION 24 [Clustering fails if $\lambda_2 \notin \mathrm{rg}(d)$ is violated]. *Assume that $\sigma(U) = \{0\} \cup \mathrm{rg}(d)$ with the eigenvalue 0 having multiplicity 1, and that the probability distribution $P$ on $\mathcal{X}$ has no point masses. Then the sequence of second eigenvalues of $\frac{1}{n}L_n$ converges to $\min_{x \in \mathcal{X}} d(x)$. The corresponding eigenfunction will approximate the characteristic function of some $x \in \mathcal{X}$ with $d(x) = \min_{x \in \mathcal{X}} d(x)$ or a linear combination of such functions.*

PROOF. It is a standard fact (Chatelin [13]) that for each $\lambda$ inside the continuous spectrum $\mathrm{rg}(d)$ of $U$ there exists a sequence of functions $(f_n)_n$ with $\|f_n\| = 1$ such that $\|(U - \lambda I)f_n\| \to 0$. Hence, for each precision $\varepsilon > 0$, there exists a function $f_\varepsilon$ such that $\|(U - \lambda I)f_\varepsilon\| < \varepsilon$. This means that for a computer with machine precision $\varepsilon$, the function $f_\varepsilon$ appears to be an eigenfunction with eigenvalue $\lambda$. Thus, with a finite precision calculation, we cannot distinguish between eigenvalues and the continuous spectrum of an operator. A similar statement is true for the eigenvalues of the empirical approximation $U_n$ of $U$. To make this precise, we consider a sequence $(f_n)_n$ as follows. For given $\lambda \in \mathrm{rg}(d)$, we choose some $x_\lambda \in \mathcal{X}$ with $d(x_\lambda) = \lambda$. Define $B_n := B(x_\lambda, \frac{1}{n})$ as the ball around $x_\lambda$ with radius $1/n$ (note that $B_n$ does not depend on the sample), and choose some $f_n \in C(\mathcal{X})$ which is constant 1 on $B_n$ and constant 0 outside $B_{n-1}$. It can be verified by straight forward arguments that this sequence has the property that for each machine precision $\varepsilon$ there exists some $N \in \mathbb{N}$ such that, for $n > N$, we have $\|(U_n - \lambda I)f_n\| \leq \varepsilon$ a.s. By Proposition 8 we can conclude that

$$\left\| \left( \frac{1}{n}L_n - \lambda I \right)(f(X_1), \ldots, f(X_n))^t \right\| \leq \varepsilon \qquad \text{a.s.}$$

Consequently, if the machine precision of the numerical eigensolver is $\varepsilon$, then this expression cannot be distinguished from 0, and the vector $(f(X_1), \ldots, f(X_n))^t$ appears to be an eigenvector of $\frac{1}{n}L_n$ with eigenvalue $\lambda$. As this construction holds for each $\lambda \in \mathrm{rg}(d)$, the smallest nonzero "eigenvalue" discovered by the eigensolver will be $\lambda_2 := \min_{x \in \mathcal{X}} d(x)$. If $x_{\lambda_2}$ is the unique point in $\mathcal{X}$ with $d(x_{\lambda_2}) = \lambda_2$, then the second eigenvector of $\frac{1}{n}L_n$ will converge to the delta-function at $x_{\lambda_2}$. If there are several points $x \in \mathcal{X}$ with $d(x) = \lambda_2$, then the "eigenspace" of $\lambda_2$ will be spanned by the delta-functions at all those points. In this case, the eigenvectors of $\frac{1}{n}L_n$ will approximate one of those delta-functions, or a linear combination thereof. $\square$

As a side remark, note that as the above construction holds for all elements $\lambda \in \mathrm{rg}(d)$, eventually the whole interval $\mathrm{rg}(d)$ will be populated by eigenvalues of $\frac{1}{n}L_n$.

So far we have seen that there exist examples where the assumption $\lambda \notin \mathrm{rg}(d)$ in Theorem 21 is violated, and that in this case the corresponding eigenfunction does not contain any useful information for clustering. This situation is aggravated by the fact that the condition $\lambda \notin \mathrm{rg}(d)$ can only be verified if the operator $U$, and hence, the probability distribution $P$ on $\mathcal{X}$, is known. As this is not the case in the standard setting of clustering, it is impossible to know whether the condition $\lambda \notin \mathrm{rg}(d)$ is true for the eigenvalues in consideration or not. Consequently, not only spectral clustering can fail in certain situations, but we are unable to check whether this is the case for a given application of clustering or not. The least thing one should do if one really wants to use unnormalized spectral clustering is to estimate the critical region $\mathrm{rg}(d)$ by $[\min_i d_i/n, \max_i d_i/n]$ and check whether the relevant eigenvalues of $\frac{1}{n}L_n$ are inside or close to this interval or not. This observation then gives an indication whether the results obtained can considered to be reliable or not.

Finally, we want to show that such problems as described above do not only occur in pathological examples, but they can come up for many similarity functions which are often used in practice.

PROPOSITION 25 (Finite discrete spectrum for analytic similarity). *Assume that $\mathcal{X}$ is a compact subset of $\mathbb{R}^n$, and the similarity function $k$ is analytic in a neighborhood of $\mathcal{X} \times \mathcal{X}$. Let $P$ be a probability distribution* on $\mathcal{X}$ which has an analytic density function. Assume that the set $\{x^* \in \mathcal{X}; d(x^*) = \min_{x \in \mathcal{X}} d(x)\}$ is finite. Then $\sigma(U)$ has only finitely many eigenvalues outside rg($d$).

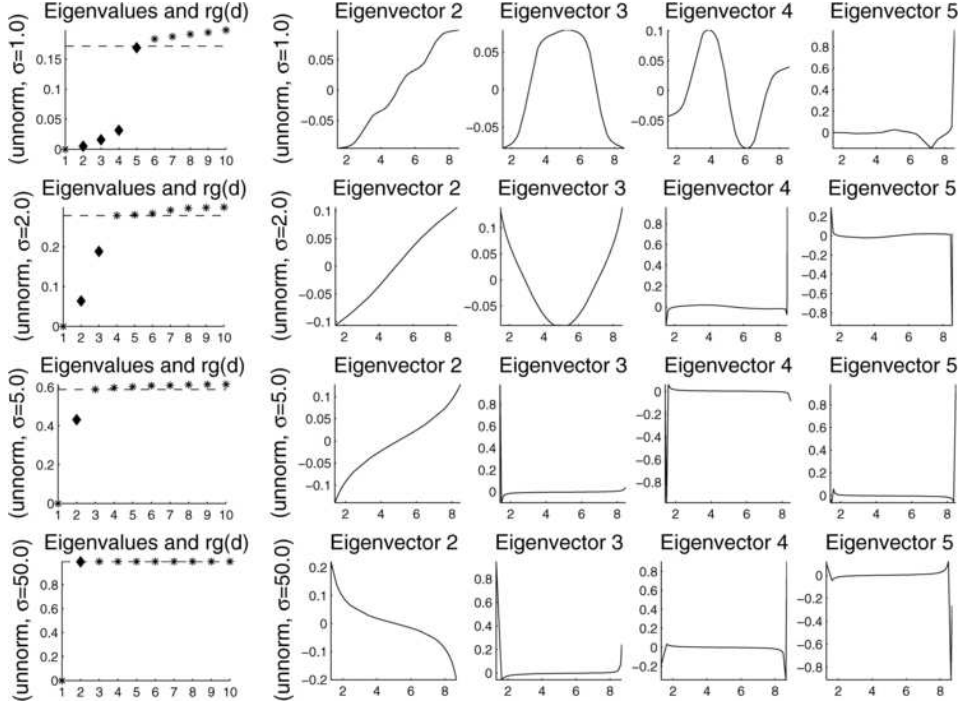

Fig. 1. *Eigenvalues and eigenvectors of the unnormalized Laplacian. Eigenvalues within* rg($d_n$) *and the trivial first eigenvalue 0 are plotted as stars, the "informative" eigenvalues below* rg($d_n$) *are plotted as diamonds. The dashed line indicates* $\min d_n(x)$. *The parameters are* $\sigma = 1$ *(first row),* $\sigma = 2$ *(second row)* $\sigma = 5$ *(third row), and* $\sigma = 50$ *(fourth row).*

This proposition is a special case of results on the discrete spectrum of the generalized Friedrichs model which can be found, for example, in Lakaev [32], Abdullaev and Lakaev [1] and Ikromov and Sharipov [26]. In those articles, the authors only consider the case where $P$ is the uniform distribution, but their proofs can be carried over to the case of analytic density functions.

8.2. *Empirical results.* To illustrate what happens for unnormalized spectral clustering if the condition $\lambda \notin$ rg($d$) is violated, we want to analyze empirical examples and compare the eigenvectors of unnormalized and normalized graph Laplacians. Our goal is to show that problems can occur in examples which are highly relevant to practical applications. As data space, we choose $\mathcal{X} = \mathbb{R}$ with a density which is a mixture of four Gaussian with means 2, 4, 6 and 8, and the same standard deviation 0.25. This density consists of four very well separated clusters, and it is so simple that every clustering algorithm should be able to identify the clusters. As similarity function we choose the Gaussian kernel function $k(x, y) = \exp(-\|x - y\|^2/\sigma^2)$, which is the similarity function most widely used in applications of spectral clustering. It is difficult to prove analytically how many eigenvalues will lie below $\mathrm{rg}(d)$; by Proposition 25, we only know that they are finitely many. However, in practice, it turns out that "finitely many" often means "very few," for example, two or three.

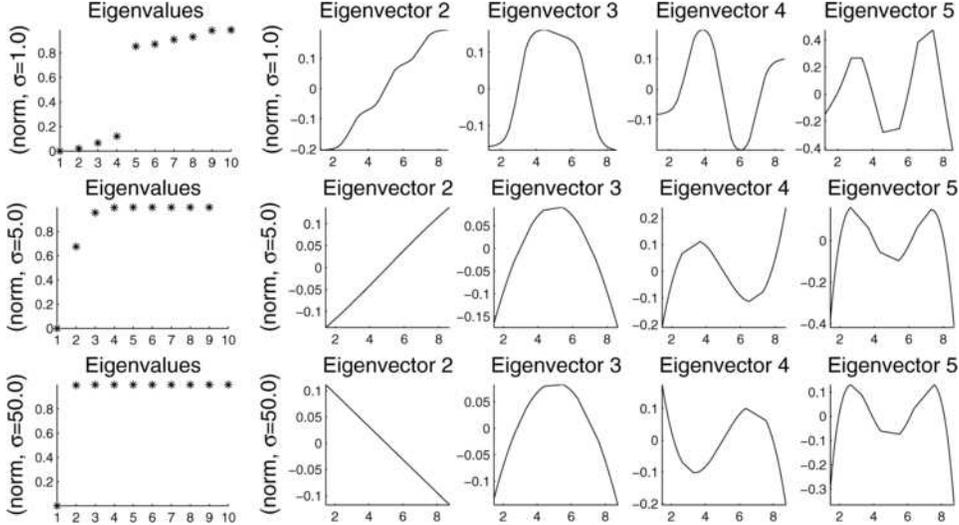

Fig. 2.  *Eigenvalues and vectors of the normalized Laplacian for* $\sigma = 1$, $\sigma = 5$ *and* $\sigma = 50$.

In Figures 1 and 2 we show the eigenvalues and eigenvectors of the normalized and unnormalized Laplacians, for different values of the kernel width parameter $\sigma$. To obtain those plots, we drew 200 data points at random from the mixture of Gaussians, computed the graph Laplacians based on the Gaussian kernel function, and computed its eigenvalues and eigenvectors. In the unnormalized case we show the eigenvalues and vectors of $L_n$, in the normalized case those of the matrix $L_n$. In each case we then plot the first 10 eigenvalues ordered by size (i.e., we plot $i$ vs. $\lambda_i$), and the eigenvectors as functions on the data space (i.e., we plot $X_i$ vs. $v_i$). In Figure 1 we show the behavior of the unnormalized graph Laplacian for various values of $\sigma$. We can observe that the larger the value of $\sigma$ is, the more the eigenvalues move toward the range of the degree function. For eigenvalues which are safely below this range, the corresponding eigenvectors are non-trivial, and thresholding them at 0 leads to a correct split between different clusters in the data (recall that the clusters are centered around 2, 4, 6 and 8). For example, in case of the plots in the first row of Figure 1, thresholding Eigenvector 2 at 0 separates the first two from the second two clusters, thresholding Eigenvector 3 separates clusters 1 and 4 from the clusters 2 and 3, and Eigenvector 4 separates clusters 1 and 3 from clusters 2 and 4. However, for eigenvalues which are very close to or inside $\mathrm{rg}(d_n)$, the corresponding eigenvector is close to a Dirac vector. In Figure 2 we show eigenvalues and eigenvectors of the normalized Laplacian. We can see that, for all values of $\sigma$, all eigenvectors are informative about the clustering, and no eigenvector has the form of a Dirac function. This is even the case for extreme values as $\sigma = 50$.

**9. Conclusion.** In this article we investigated the consistency of spectral clustering algorithms by studying the convergence of eigenvectors of the normalized and unnormalized Laplacian matrices on random samples. We proved that, under standard assumptions, the first eigenvectors of the normalized Laplacian converges to eigenfunctions of some limit operator. In the unnormalized case, the same is only true if the eigenvalues of the limit operator satisfy certain properties, namely, if these eigenvalues lie below the continuous part of the spectrum. We showed that in many examples this condition is not satisfied. In those cases, the information provided by the corresponding eigenvector is misleading and cannot be used for clustering.

This leads to two main practical conclusions about spectral clustering. First, from a statistical point of view, it is clear that normalized rather than unnormalized spectral clustering should be used whenever possible. Second, if for some reason one wants to use unnormalized spectral clustering, one should try to check whether the eigenvalues corresponding to the eigenvectors used by the algorithm lie significantly below the continuous part of the spectrum. If that is not the case, those eigenvectors need to be discarded, as they do not provide information about the clustering.

U. VON LUXBURG
MAX PLANCK INSTITUTE
  FOR BIOLOGICAL CYBERNETICS
SPEMANNSTR. 38
72076 TÜBINGEN
GERMANY
E-MAIL: ulrike.luxburg@tuebingen.mpg.de

M. BELKIN
DEPARTMENT OF COMPUTER SCIENCE
  AND ENGINEERING
OHIO STATE UNIVERSITY
2015 NEIL AVENUE
COLUMBUS, OHIO 43210
USA
E-MAIL: mbelkin@cse.ohio-state.edu

O. BOUSQUET
PERTINENCE
32 RUE DES JEÛNEURS
F-75002 PARIS
FRANCE
E-MAIL: olivier.bousquet@pertinence.com